\documentclass[hidelinks]{siamart171218}
\usepackage{smalljeff}
\usepackage{combelow}
\pgfplotsset{compat=newest}
\usepgfplotslibrary{groupplots}

\newcommand{\Htwo}{{\mathcal{H}_2}}

\def\stripzero#1{\expandafter\stripzerohelp#1}
\def\stripzerohelp#1{\ifx 0#1\expandafter\stripzerohelp\else#1\fi}

\newcommand{\TheTitle}{Least Squares Rational Approximation}
\newcommand{\TheAuthors}{Jeffrey M. Hokanson and Caleb C. Magruder}
\headers{Least Squares Rational Approximation}{\TheAuthors}

\title{{\TheTitle}\thanks{Submitted to the editors 15 June 2018.
\funding{The first author's work is partially supported by Department of Defense, 
Defense Advanced Research Project Agency's program Enabling Quantification of Uncertainty in Physical Systems. 
}}}

\author{Jeffrey M. Hokanson\thanks{
	Department of Computer Science, 
	University of Colorado Boulder,
	1111 Engineering Dr, Boulder, CO 80309,
	(\email{Jeffrey.Hokanson@colorado.edu}).}
	\and Caleb C. Magruder\thanks{
	Department of Mathematics,
	Tufts University,
	503 Boston Avenue, 
	Bromfield-Pearson, 
	Medford, MA 02155
	(\email{Caleb.Magruder@tufts.edu})}
}

\begin{document}
\maketitle 
\begin{abstract}
Rational approximation appears in many contexts throughout science and engineering,
playing a central role in linear systems theory, special function approximation,
and many others.
There are many existing methods for solving the rational approximation problem,
from fixed point methods like the Sanathanan-Koerner iteration and Vector Fitting,
to partial interpolation methods like Adaptive Anderson Antoulas (AAA).
While these methods can often find rational approximations with a small residual norm, 
they are unable to find optimizers with respect to a weighted $\ell_2$ norm with a 
square dense weighting matrix. 
Here we develop a nonlinear least squares approach 
constructing rational approximations with respect to this norm.
We explore this approach using two parameterizations of rational functions:
a ratio of two polynomials and a partial fraction expansion.
In both cases, we show how we can use Variable Projection (VARPRO)
to reduce the dimension of the optimization problem.
As many applications seek a real rational approximation
that can be described as a ratio of two real polynomials,
we show how this constraint can be enforced in both parameterizations.
Although this nonlinear least squares approach often converge to suboptimal local minimizers,
we find this can be largely mitigated by initializing the algorithm using 
the poles of the AAA algorithm applied to the same data.
This combination of initialization and nonlinear least squares enables
us to construct rational approximants using dense and potentially ill-conditioned weight matrices
such as those that appear as a step in new $\Htwo$ model reduction algorithm
recently developed by the authors. 
\end{abstract}

\begin{keywords}
	rational approximation,
	nonlinear least squares,
	variable projection
\end{keywords}
\begin{AMS}	
	26C15, 
	41A20, 
	90C53 
\end{AMS}
\begin{DOI}

\end{DOI}

\section{Introduction}
Rational approximation plays a role in several applications in science and engineering;
for example, rational approximations are a critical component in $\Htwo$-model reduction~\cite{DGB15},
can be used in special function computation~\cite[sec.~9.2]{GST07},
and many others~\cite[Chap.~23]{Tre13}. 
Generally posed, the goal of rational approximation is to mimic a function $f:\C\to \C$
by a degree $(m,n)$ rational function $r(z): \C\to \C$ 
\begin{equation}
	f \approx r \in \set R_{m,n}(\C) := 
	\left\lbrace \frac{p}{q} : p \in \set P_m(\C), \ q \in \set P_n(\C) \backslash \{0\} 
	\right\rbrace
\end{equation}
where $\set P_m(\C)$ denotes the set of polynomials of degree $m$ with coefficients in $\C$.
There are several senses in which we might seek to construct a rational approximation.
For example, Pad\'e approximation~\cite{BG96} chooses $r$ to match the first $m+n$ derivatives of $f$
at some point $\widehat z\in \C$:
\begin{equation}
	r \in \set R_{m,n}(\C) \quad \text{such that} \quad
	f^{(k)} (\widehat z) = r^{(k)} (\widehat z) \quad \forall k = 0, 1,\ldots, m+n.
\end{equation}
In special function approximation, 
the goal is often to construct a minimax rational approximation~\cite{EW76}
that minimizes the maximum mismatch over a set $\set Z \subset \C$:
\begin{equation}\label{eq:ratfit}
	\minimize_{r \in \set R_{m,n}(\C) } \ \sup_{z \in \set Z} \left| f(z) -  r(z) \right|.
\end{equation}
In this paper we seek to construct a \emph{least squares rational approximation}
over a discrete set of $N$ points $\set Z\subset \C$
in a weighted $\ell_2$ norm with a dense weight matrix $\ma W \in \C^{N\times N}$:
\begin{equation}\label{eq:ls_ratfit}
	\minimize_{r \in \set R_{m,n}(\C) }
		\| \ma W[ f(\set Z) - r(\set Z) ]\|_2^2
	\quad \text{where} \quad  
		f(\set Z) := \begin{bmatrix} f(z_1) \\ \vdots \\ f(z_N) \end{bmatrix}, \
		r(\set Z) := \begin{bmatrix} r(z_1) \\ \vdots \\ r(z_N) \end{bmatrix}.
\end{equation}
Our motivation for studying this weighted least squares rational approximation 
comes from a new $\Htwo$ model reduction algorithm developed by the authors~\cite{HM18x}
where a rational approximation of this form appears at each step of the algorithm
with a weight matrix $\ma W$ that is the inverse matrix square root of a Cauchy matrix.
Existing algorithms for rational approximation cannot incorporate this non-diagonal 
weight matrix, leading us to develop an algorithm to solve~\cref{eq:ls_ratfit}
based standard nonlinear least squares techniques.

There are a variety of existing algorithms for rational approximation.
For example as discussed in \cref{sec:background:partial},
the Loewner\footnote{also spelled L\"owner} framework of Anderson and Antoulas~\cite{AA86}
for rational interpolation has been extended
by Nakatsukasa, S\`ete, and Trefethen~\cite{NST18x} to rational approximation problem
in the Adaptive Anderson-Antoulas (AAA) algorithm;
however this does not minimize the nonlinear least squares problem~\cref{eq:ls_ratfit}
and does not incorporate a weighting matrix.
Similarly, as discussed in \cref{sec:background:fixed},
there are fixed point methods such as 
the Sanathanan-Koerner (SK) iteration~\cite{SK63} and 
the \emph{Vector Fitting} algorithm of Gustavsen and Semlyen~\cite{GS99}
which have fixed points nearby minimizers of the nonlinear least squares problem~\cref{eq:ls_ratfit}.
These methods can incorporate a diagonal weighting matrix, 
but the dense weighting matrix $\ma W$ required for the $\Htwo$ model reduction problem.

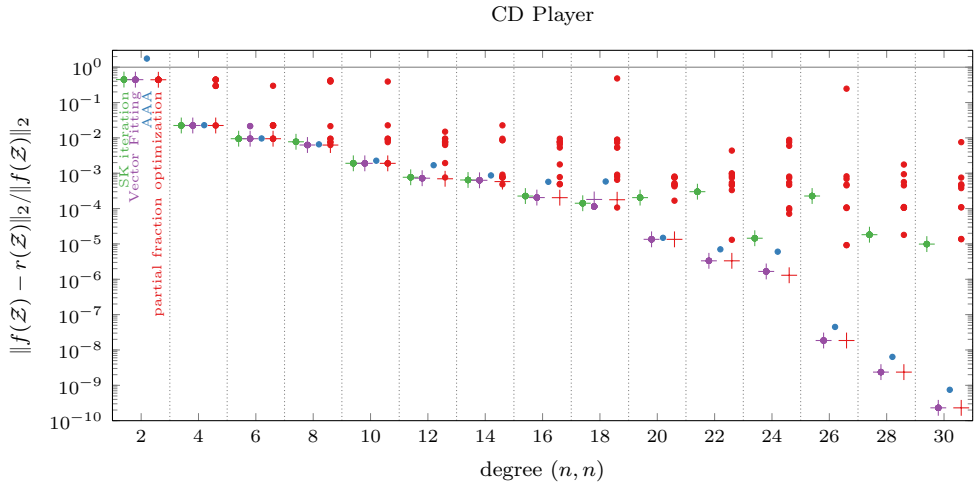
\begin{figure}

\begin{tikzpicture}
\begin{axis}[
	width = \textwidth,
	height = 0.5\textwidth,
	ymode = log,
	xmin = 1,
	xmax = 31,
	ymax = 3e0,
	ymin = 1e-10,
	ytickten = {-10,-9,...,0},
	xlabel = {degree $(n,n)$},
	ylabel = {$\| f(\set Z) - r(\set Z)\|_2/\| f(\set Z)\|_2$},
	title = {CD Player},	
]
	\foreach \rn in {02, 04, 06, 08, 10, 12, 14, 16, 18, 20, 22, 24, 26, 28, 30}{
		\addplot[colorbrewerA3, mark=*, only marks, mark size = 1] table
			 [x expr = {\stripzero{\rn} - 0.6}, y = res_sk ]
			{data/fig_local_minima_\rn.dat};

		\addplot[colorbrewerA4, mark=*, only marks, mark size = 1] table
			 [x expr = \stripzero{\rn} - 0.2, y = res]
			{data/fig_local_minima_vectfit3_\rn.dat};
		

		\addplot[colorbrewerA1, mark=*, only marks, mark size = 1] table 
			[x expr = \stripzero{\rn}+0.6, y = res_pr ]
			{data/fig_local_minima_\rn.dat};
	}
	\foreach \rn in {1,3,...,31}{
		\addplot[gray, densely dotted] coordinates {(\rn,1e-10) (\rn, 3)};	
	}

	\addplot[colorbrewerA3, mark = +, only marks, mark size = 3, ] table
		[x expr = \thisrow{n} - 0.6, y = res_sk]
		{data/fig_local_minima.dat};

	\addplot[colorbrewerA4, mark = +, only marks, mark size = 3, ] table
		[x expr = \thisrow{n} - 0.2, y = res]
		{data/fig_local_minima_vectfit3_aaa_init.dat};
	
	\addplot[colorbrewerA1, mark = +, only marks, mark size = 3, ] table
		[x expr = \thisrow{n} + 0.6, y = res_pr]
		{data/fig_local_minima.dat};
	
	\addplot[colorbrewerA2, mark = *, only marks, mark size = 1, ] table
		[x expr = \thisrow{n} + 0.2, y = res_aaa]
		{data/fig_local_minima.dat};
	
	\addplot[colorbrewerA2, mark = +, only marks, mark size = 1, ] table
		[x expr = \thisrow{n} + 0.2, y = res_aaa]
		{data/fig_local_minima.dat};

	\draw [gray] (0,1) -- (32,1);
	\draw [colorbrewerA3] (1.4, 0.4) node [anchor = east, rotate = 90] {\tiny SK iteration};
	\draw [colorbrewerA4] (1.8, 0.4) node [anchor = east, rotate = 90] {\tiny Vector Fitting};
	\draw [colorbrewerA2] (2.2, 0.4) node [anchor = east, rotate = 90] {\tiny AAA};
	\draw [colorbrewerA1] (2.6, 0.4) node [anchor = east, rotate = 90] {\tiny partial fraction optimization};
\end{axis}
\end{tikzpicture}

\caption[Illustration of Local Minima]{
	Spurious local minima are a significant problem 
	when building a rational approximation via optimization
	compared to the SK iteration and Vector Fitting described in \cref{sec:background:fixed}.
	Here each dot \tikz{\draw[black] plot[mark=*, mark size = 1](0,0);}
	shows the normalized residual of rational approximants generated from
	ten different initializations of each algorithm.
	Both the SK iteration and Vector Fitting frequently converge to a rational approximant 
	with similar residual norm, whereas our optimization approach strongly depends on the starting
	rational function.
	However, by initializing using AAA,
	denoted by crosses \tikz{\draw[black] plot[mark=+, mark size = 3] (0,0);},
	we are able to mitigate this dependence on initial condition
	and find a local minimizer that is comparable to that of Vector Fitting.
	The performance of the SK iteration after $n\ge 18$ is caused by ill-conditioning as illustrated in 
	\cref{fig:fixed_point}.	
	In this example $f$ is the $(1,1)$ entry of the transfer function of the CD player model~\cite{CD02}
	evaluated at 1000 points evenly sampled on the imaginary axis between $-10^3i$ and $10^3i$.
}
\label{fig:local_minima}
\end{figure}

One might ask: why not use nonlinear least squares methods?
To use this approach, we first need to specify a parameterization 
for the rational approximant $r$.
Although there are many potential parameterizations, here we focus on two:
a polynomial parameterization and a partial fraction parameterization.
In the polynomial parameterization,
we define $r$ by the coefficients $\ve a\in \C^m$ and $\ve b\in \C^n$ 
of the numerator and denominator polynomials
expressed in bases $\lbrace \phi_k\rbrace_{k=0}^m \subset \set P_m(\C)$
and $\lbrace \psi_k \rbrace_{k=0}^n \subset \set P_n(\C)$
\begin{equation}\label{eq:polybasis}
	r(z; \ve a, \ve b) := 
		\frac{p(z; \ve a)}{q(z; \ve b)} = 
	\frac{\sum_{k=0}^{m} a_k \phi_k(z)}{\sum_{k=0}^n b_k \psi_k(z)}.
\end{equation}
In a partial fraction parameterization, 
which is limited to degree $(m,n)$ rational functions where $m\ge n-1$,
we define $r$ as a sum of degree $(0,1)$ rational functions
described by their poles $\ve\lambda \in \C^n$ and residues $\ve \rho\in \C^n$
plus an additional set of polynomial coefficients $\ve c \in \C^{m-n}$
\begin{equation}\label{eq:polebasis}
	r(z; \ve \lambda, \ve \rho, \ve c) := 
		\sum_{k=1}^n \frac{\rho_k}{z-\lambda_k} + \sum_{k=0}^{m-n} c_k \varphi_k(z).
\end{equation}
Sections \ref{sec:polyopt} and \ref{sec:propt} provide formulas for the residual and Jacobian 
for these two parameterizations.
However, with either parameterization the challenge with this approach is spurious local minima.
As illustrated in \cref{fig:local_minima} when optimizing in the partial fraction parameterization
starting from different initialization, the algorithm finds different local minimizers.
Moreover these minimizers frequently have a larger residual norm than the solutions generated by 
the Vector Fitting and the SK iteration before $n > 14$ when numerical instability emerges. 
This explains the relative infrequency with which rational approximation is treated using an optimization approach;
we are only aware of one paper by Lefteriu and Antoulas where this approach is briefly described~\cite{LA13}.
One approach to mitigate the issue of spurious local minima is to find an effective initialization.
Here we advocate using the AAA algorithm as an initialization for for a standard Gauss-Newton method
with a backtracking line search (see, e.g.,~\cite[sec.~10.3]{NW06};
we use this combination to construct the remainder of our examples.
As evidenced in \cref{fig:local_minima},
coupling this initialization approach with Gauss-Newton yields better minimizers than random initialization,
and the residual norm associated with these optimizers is comparable to that generated by Vector Fitting.

Additionally, in this paper we also address how to construct a real rational approximation,
where $r$ is a real rational function
\begin{equation}
	r \in \set R_{m,n}(\R) := 
	\left\lbrace \frac{p}{q} : p \in \set P_m(\R), \ q \in \set P_n(\R) \backslash \{0\} 
	\right\rbrace
\end{equation}
where $\set P_n(\R)$ is the space of degree $n$ polynomials with real coefficients.
This constraint is frequently present in model reduction context
since the transfer function $f(z) = \ve c^* (\ma I z - \ma A)^{-1}\ve b$
is a real rational function if $\ma A$, $\ve b$, and $\ve c$ are real.
Although imposing this constraint in the polynomial parameterization is straightforward,
doing so in the partial fraction expansion requires more care
which we discuss in \cref{sec:propt:real}.
 
In the remainder of this paper we first 
describe AAA algorithm, the SK iteration, and Vector Fitting in \cref{sec:background}
and discuss their numerical proprieties.
Then in sections \ref{sec:polyopt} and \ref{sec:propt} we derive
the residual and Jacobian for the polynomial and partial fraction parameterizations
using Variable Projection (VARPRO)~\cite{GP73} to pose the optimization problem only over the nonlinear parameters
and also showing how to enforce the constraint that the rational approximant is real.
We conclude that due to poor conditioning, the partial fraction parameterization is preferable
unless the goal is to build a rational approximation of degree $(m,n)$ where $m< n-1$,
which can not be expressed in this parameterization.
Then in \cref{sec:examples} we provide an example of the effect of imposing the real constraint
and evaluate the performance of our algorithm in an example mimicking a step of the projected $\Htwo$ algorithm.
Finally we discuss extending our optimization approach to vector and matrix valued output data in 
\cref{sec:discussion}.

\section{Existing Methods for Rational Approximation\label{sec:background}}
In this section we describe two popular classes of algorithms for rational approximation:
those based on the Loewner framework originating in the work of Antoulas and Anderson~\cite{AA86}, 
such as the Adaptive Antoulas Anderson (AAA) algorithm~\cite{NST18x},
and fixed point iterations, such as the Sanathanan-Koerner iteration~\cite{SK63}
and Vector Fitting~\cite{GS99}.
While these methods are successful in consistently obtaining a rational approximation
with a small residual as illustrated in \cref{fig:local_minima},
none of these methods can incorporate a dense square mass matrix.
Moreover, our numerical experiments suggest that
that the rational approximations that these methods generate 
do not satisfy the first order necessary conditions
for the least squares rational approximation problem~\cref{eq:ratfit}.
In the remainder of this section we will briefly derive each method, 
illustrating that each method uses a similar trick---%
multiplying by the denominator of the polynomial, effectively `linearizing' the problem---%
and discuss how this affects the ability of the algorithm to obtain a least squares estimate.
Although these methods do not provide least squares estimates,
they are capable of providing rational interpolants (in exact arithmetic) when both $r$ and $f$
are degree $(m,n)$ rational functions.

\subsection{Loewner Framework\label{sec:background:partial}}
The original work by Anderson and Antoulas 
presented a technique for determining if a rational interpolant of a specified degree exists
and, if so, gave a formula for such a rational interpolant~\cite[eq.~(2.11)]{AA86}.
The central feature of this analysis is a (generalized) Loewner matrix,
defined through the input-output pairs $(z_j, f(z_j))$.
Although the original derivation permitted interpolation including arbitrary orders of derivatives,
here we will describe a simplification for the rational approximation problem
called the Adaptive Antoulas Anderson (AAA) algorithm developed by Nakatsukasa, S\`ete, and Trefethen~\cite{NST18x}.

A key feature of methods in the Loewner framework is splitting 
the sample points $\set Z$ into two disjoint sets:
$\hset Z := \lbrace \widehat{z}_k\rbrace_{k=0}^n \subset \set Z$
and $\cset Z :=\lbrace \widecheck{z}_k\rbrace_{k=0}^{N-n} \subset \set Z$
where $\hset Z \cup \cset Z = \set Z$.
Then, the key step this approach takes to relax the rational approximation problem is
forcing the rational approximant $r$ to interpolate at the values in $\hset Z$.
This results in a simple method to find the rational partial-interpolant---%
namely, the singular value decomposition (SVD).

To see this in the context of the AAA algorithm, 
we express the rational approximation $r$ using the same Lagrange basis for numerator and denominator,
restricting this approach to building degree $(n,n)$ rational approximants.
In particular, we will use a Lagrange basis expressed in an unweighted barycentric form~\cite[eq.~(3.3)]{BT04}
with Lagrange nodes $\lbrace \widehat{z}_j\rbrace_{j=0}^n=\set Z$ where we will later force interpolation:
\begin{equation}
	\phi_k(z) = (z - \widehat z_k)^{-1} \ell(z),
	\qquad
	\ell(z) = \prod_{k=0}^n (z - \widehat z_k)
\end{equation}
Then, writing the rational approximant $r$ (cf.~\cref{eq:polybasis})
allows us to cancel the common factor $\ell(z)$:
\begin{equation}
	r(z; \ve a, \ve b) =
		\frac{ \ell(z)\sum_{k=0}^n a_k (z-\widehat z_k)^{-1}}
		{\ell(z)\sum_{k=0}^n b_k  (z-\widehat z_k)^{-1}} 
		=\frac{ \sum_{k=0}^n a_k (z-\widehat z_k)^{-1}}
		{\sum_{k=0}^n b_k (z-\widehat z_k)^{-1}}
\end{equation}
Next, invoking the (suboptimal) assumption that $r$ interpolates $f$ on $\hset Z$, 
we require that (after removing a removable singularity),
\begin{equation}\label{eq:aaa_interp}
	r(\widehat z_k; \ve a, \ve b) = \frac{a_k}{b_k} = f(\widehat z_k).
\end{equation}
Hence if we assume $b_k\ne 0$, we can set $a_k := b_k f(\widehat z_k)$
yielding an expression for our rational approximant only in terms of $\ve b$:
\begin{equation}
	r(z; \ve b) := 
		\frac{ \sum_{k=0}^n f(\widehat z_k) b_k (z-\widehat z_k)^{-1}}
		{\sum_{k=0}^n b_k (z-\widehat z_k)^{-1}}.
\end{equation}
This is related to the second form of the barycentric formula~\cite[eq.~(4.2)]{BT04},
where in the case of polynomial approximation $b_k$ is fixed $b_k = \prod_{k\ne j}(\widehat z_k - \widehat z_j)$;
this expression is also called the \emph{rational barycentric formula}~\cite[eq.~(1.7)]{Ion13}.

At this point we still need to find a choice of $\ve b$ 
such that $r$ approximates $f$ well on the remainder of the points in $\set Z$, namely $\cset Z$.
Ideally, we would solve the nonlinear least squares problem:
\begin{equation}\label{eq:aaa_ideal}
	\minimize_{\ve b, \|\ve b\|_2=1} \| f(\cset Z) - r(\cset Z; \ve b)\|_2,
\end{equation}
however, this is still a challenging nonlinear least squares problem.
Instead, if we multiply through by the denominator,
introducing a second modification of the optimization problem,
we find that $\ve b$ now appears linearly in each row
\begin{multline}\label{eq:aaa_linearize}
	f(\widecheck z_j) - \frac{ \sum_{k=0}^n f(\widehat z_k) b_k (\widecheck z_j - \widehat z_k)^{-1}}
		{\sum_{k=0}^n b_k (\widecheck z_j - \widehat z_k)^{-1}}
	\Rightarrow
	f(\widecheck z_j) \sum_{k=0}^n \frac{b_k}{\widecheck z_j - \widehat z_k}
	- \sum_{k=0}^n \frac{b_k f(\widehat z_k)}{\widecheck z_j - \widehat z_k}.
\end{multline}
After this modification,
we find $\ve b$ as the smallest singular value of the Loewner matrix $\ma L\in \C^{(N-n-1) \times (n+1)}$
built from the input output pairs:
\begin{equation}\label{eq:aaa}
	\minimize_{\ve b, \|\ve b\|_2 = 1} \| \ma L \ve b\|_2,
	\quad \ma L = \begin{bmatrix}
		\frac{f(\widecheck z_1) - f(\widehat z_0)}{\widecheck z_1 - \widehat z_0} &
		\ldots &
		\frac{f(\widecheck z_1) - f(\widehat z_n)}{\widecheck z_1 - \widehat z_n} \\ 
		\vdots & & \vdots \\
		\frac{f(\widecheck z_{N-n-1}) - f(\widehat z_0)}{\widecheck z_{N-n-1} - \widehat z_0} &
		\ldots &
		\frac{f(\widecheck z_{N-n-1}) - f(\widehat z_n)}{\widecheck z_{N-n-1} - \widehat z_n} \\ 
	\end{bmatrix}.
\end{equation}
The net result of these approximations is an easy approach for finding a rational approximant,
but one that is necessary suboptimal with respect to the $\ell_2$ norm
due to the interpolation condition~\cref{eq:aaa_interp}
and multiplication by the denominator in~\cref{eq:aaa_linearize}.
However, as illustrated in \cref{fig:aaa},
this approach will yield increasingly good rational approximants
as measured in residual norm with increasing degree,
but ones that are outperformed by our optimization based approach for building rational approximants. 

\begin{figure}
\begin{tikzpicture}
\begin{groupplot}[
	group style = {group size = 2 by 1, horizontal sep = 3.5em}, 
	width=0.5\textwidth,
	height = 0.35\textwidth,
	ymode = log, 
	xmin = 0, 
	xmax = 30.5,
	xtick = {1,5,10,15,...,50},
	xlabel = {degree $(m,m)$},
	ytickten = {-10,-8,...,2},
	]
	\nextgroupplot[ymin = 1e-10, ymax = 1e1, ylabel = normalized residual, title = CD Player]
	\addplot[colorbrewerA1, mark=*, only marks, mark size=1] table 
		[x = m, y = res_pr]
		{data/fig_aaa_cdplayer.dat};
	\addplot[colorbrewerA2, mark=o, only marks, mark size=1.5] table 
		[x = m, y = res_aaa]
		{data/fig_aaa_cdplayer_matlab.dat};
	
	\nextgroupplot[ymin = 1e-10, ymax = 1e1, title = $\tan(256z)$, xmax=50.5]
	\addplot[colorbrewerA1, mark=*, only marks, mark size=1] table 
		[x = m, y = res_pr]
		{data/fig_aaa_tan256.dat};
	\addplot[colorbrewerA2, mark=o, only marks, mark size=1.5] table 
		[x = m, y = res_aaa]
		{data/fig_aaa_tan256_matlab.dat};
\end{groupplot}
\end{tikzpicture}
\caption[AAA Comparision]{
Comparison of the residual $\|f(\set Z) - r(\set Z)\|_2$ 
where the rational approximation $r$ is generated by
AAA (\tikz{\draw[colorbrewerA2] plot[mark=o, mark size=1.5] (0,0);})
and our optimization approach using a pole residue parameterization described in \cref{sec:propt}
(\tikz{\draw[colorbrewerA1] plot[mark=*, mark size=1] (0,0);}).
The left takes $f$ from the CD Player model as described in \cref{fig:local_minima}.
The right takes $f(z) = \tan(256z)$ evaluated at 1000 points uniformly on the unit circle
following \cite[Fig.~6.4]{NST18x}.
}
\label{fig:aaa}
\end{figure}

One of the important contributions of the AAA algorithm
was providing a greedy heuristic for selecting interpolation points,
given in \cref{alg:aaa}.
The authors also discuss removing Froissart doublets---%
poles with either small residues or pole-zero pairs that nearly cancel---%
an important consideration for the quality of approximation
when the norm of the residual becomes small.
\begin{algorithm}[t]
\begin{minipage}{\linewidth}
\begin{algorithm2e}[H]
\Input{Input output pairs $\lbrace (z_j, f(z_j))\rbrace_{j=1}^N$, desired degree approximant $(n,n)$}
\Output{Rational approximation $r(z) = 
			\frac{\sum_{k=0}^\ell f(\widehat z_k) b_k(z-z_j)^{-1}}{\sum_{k=0}^\ell b_k(z - \widehat z_k)^{-1} }$}
Set residual $r_j \leftarrow f(z_j)\quad j=1,\ldots,N$\;
$\cset Z \leftarrow \set Z$\;
\For{$\ell=0,1,\ldots, n$}{
	$j\leftarrow \argmax_{j} |r_j|$\;
	Remove $z_j$ from $\cset Z$ and place in $\hset Z$\;
	Construct the Loewner matrix 
	$[\ma L]_{j,k} \leftarrow (f(\widecheck z_j) - f(\widehat z_k))/(\widecheck z_j - \widehat z_k)$\;
	Compute SVD: $\ma U \ma \Sigma \ma V^*\leftarrow \ma L$\;
	$\ve b \leftarrow [\ma V]_{\cdot,\ell}$\;
	Compute residual 
	$r_j \leftarrow \begin{cases}
		f(z_j) -\frac{\sum_{k=0}^\ell f(\widehat z_k) b_k (z_j - \widehat z_k)^{-1}}{\sum_{k=0}^\ell b_k(z_j -\widehat z_k)^{-1}}, & z_j \in \cset Z\\
		0, & z_j \in \hset Z.
		\end{cases}
	$\; 
}
\end{algorithm2e}
\vspace{-1em}
\end{minipage}
\caption{Adaptive Anderson Antoulas (AAA)}
\label{alg:aaa}
\end{algorithm}

\subsection{Fixed Point Iterations\label{sec:background:fixed}}
An alternative to the Loewner framework are fixed point iterations,
such as the \emph{Sanathanan-Koerner (SK) iteration}~\cite{SK63}
and \emph{Vector Fitting}~\cite{GS99, Gus06}.
Both iterations exploit the same trick of multiplying through by the denominator 
which was seen in~\cref{eq:aaa_linearize}
and, unlike AAA, do not require interpolation at a set of $n+1$ points.
Although the SK iteration and Vector Fitting were developed independently,
their similar underpinning has been previously discussed by Hendrickx and Dhaene~\cite{HD06}.
Here we focus on the numerical features of each algorithm,
noting the SK iteration can become ill-conditioned even when using a
polynomial basis that is well conditions, such as a Legendre basis.
However, Vector Fitting avoids this fault by working asymptotically in a partial fraction expansion.
Although both algorithms can provide good rational approximations,
neither satisfy the first order optimality criteria for least squares rational approximation
when data $f(z_i)$ is generated by a function $f$ not in $\set R_{m,n}$~\cite{Shi16}.

\subsubsection{The Sanathanan-Koerner Iteration}
Suppose we are given a basis for the numerator and denominator,
$\lbrace \phi_k\rbrace_{k=0}^m$ and $\lbrace \psi_k\rbrace_{k=0}^n$
and construct the rational approximation $r$ as in \cref{eq:polybasis}:
\begin{equation}
	r(z;\ve a, \ve b) = \frac{p(z; \ve a)}{q(z;\ve b)}= 
		\frac{\sum_{k=0}^m a_k \phi_k(z)}{\sum_{k=0}^n b_k \psi_k(z)}.
\end{equation}
If we define Vandermonde matrices $\ma \Phi \in \C^{N\times m}$ and $\ma \Psi \in \C^{N\times n}$,
we can write this optimization problem as 
\begin{equation}
	\minimize_{\ve a \in \C^{m+1}, \ve b\in \C^{n+1}} \| \ve f - \diag(\ma \Psi \ve b)^{-1} \ma \Phi \ve a\|_2,
	\quad
	[\ma \Phi]_{j,k} := \phi_k(z_j),
	\quad
	[\ma \Psi]_{j,k} := \psi_k(z_j),
\end{equation}
and $\ve f = [f(z_1),\ldots, f(z_N)]^\trans$.
One common approach to building a rational approximation prior to Sanathanan and Koerner's 1963
paper was to multiply through by the denominator, as in AAA, yielding a linear least squares problem:
\begin{equation}
	\minimize_{\ve a \in \C^{m+1}, \ve b\in \C^{n+1}} \| \diag(\ve f) \ma \Psi \ve b - \ma \Phi \ve a\|_2.
\end{equation}
The key insight of Sanathanan and Koerner was to introduce a weighting 
to correct for the wrong norm introduced by multiplying through by the denominator.
If at step $\ell$, we have coefficients $\ve a^{(\ell)}$ and $\ve b^{(\ell)}$,
the next step is chosen by solving a problem weighted by the previous denominator:
\begin{equation}\label{eq:sk_step}
	\ve a^{(\ell+1)}, \ve b^{(\ell+1)} \leftarrow
	\minimize_{\ve a \in \C^{m+1}, \ve b\in \C^{n+1}} 
		\| \diag(\ma \Psi \ve b^{(\ell)})^{-1}[\diag(\ve f) \ma \Psi \ve b - \ma \Phi \ve a]\|_2.
\end{equation} 
Then if $\ve a^{(\ell)}\to \ve a^{(*)}$ and $\ve b^{(\ell)}\to \ve b^{(*)}$,
then $\ve a^{(*)}$ and $\ve b^{(*)}$ appear to provide a least squares solution.

There is one additional choice that is left to be made:
how to fix the free scaling shared between $\ve a$ and $\ve b$.
In the original paper Sanathanan and Koerner, working in the monomial basis, 
pick the constant term to set to one that $\ve f$ is in the right hand side of the least squares problem.
Here we follow a similar approach when $\lbrace \psi_k \rbrace_{k=0}^n$ 
is an orthogonal basis of increasing degree, such a Legendre polynomials,
which in this case yields the step:
\begin{equation}\label{eq:skiter}
	\ve a^{(\ell+1)}, \ve b^{(\ell+1)} \leftarrow
	\minimize_{\ve a \in \C^{m+1}, \ve b\in \C^n} 
		\left\| \diag(\ma \Psi \ve b^{(\ell)})^{-1}
			\left[\ve f +  [\ma \Psi]_{\cdot,1:n} [\ve b]_{1:n} - \diag(\ve f)\ma \Phi \ve a\right]
		\right\|_2;
	\quad
\end{equation} 
where $b_0^{(\ell+1)} = \psi_0(0)$.
This approach is used in our implementation, given in~\cref{alg:sk}. 
Our experiments suggest that changing the normalization 
changes the fixed points of this algorithm.
This particular constraint yielded the best fixed points in terms of residual norm
of those we experimented with.

\begin{algorithm}[t]
\begin{minipage}{\linewidth}
\begin{algorithm2e}[H]
	\Input{Input-outout pairs $\lbrace (z_j, f(z_j))\rbrace_{j=1}^N$; 
			polynomial bases $\lbrace \phi_k\rbrace_{k=0}^m \in \mathcal{P}_m$
			and $\lbrace \psi_k \rbrace_{k=0}^n \in \mathcal{P}_n$ where $\psi_0$ is constant.}
	\Output{Rational approximation $r(z) = (\sum_{k=0}^m a_k \phi_k(z) )/(\sum_{k=0}^n b_k \psi_k(z) )$}
	$\ve b^{(0)} = \phi_0(0)\ve e_0$ \;
	\For{$\ell=1,2,\ldots$}{
		$\displaystyle 
			\begin{bmatrix} \ve a^{(\ell)} \\ \hve b^{(\ell)} \end{bmatrix}
			\leftarrow \minimize_{\ve x \in \C^{m+n+1} } \left\| 
			\diag(\ma \Psi \ve b^{(\ell-1)})^{-1}
			\left[\ve f +  
			\begin{bmatrix}
				-\ma \Phi & \diag(\ve f)[\ma \Psi]_{\cdot,1:n}	
			\end{bmatrix}\ve x
			\right]
			\right\|_2$\;
		$\ve b^{(\ell)} \leftarrow \begin{bmatrix} \phi_0(0) & \hve b^\trans\end{bmatrix}^\trans$\;
		\lIf{$\| \ve b^{(\ell)} - \ve b^{(\ell-1)} \|_2<$ tol}{ {\bf break} }
	}
\end{algorithm2e}
\vspace{-1em}
\end{minipage}
\caption{Generalized Sanathanan-Koerner Iteration}
\label{alg:sk}
\end{algorithm}

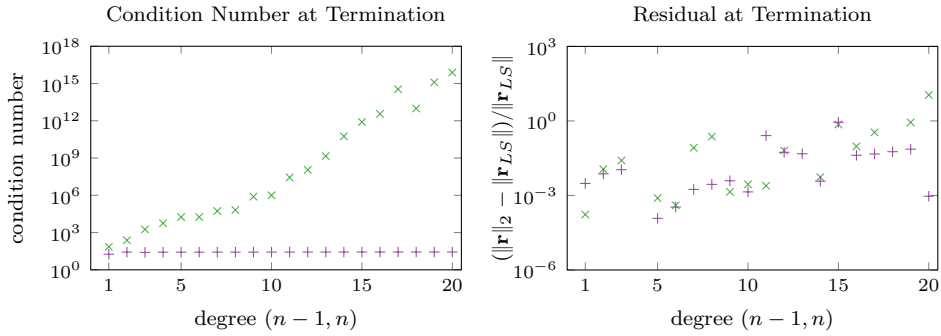
\begin{figure}
\centering
\begin{tikzpicture}
\pgfplotstableread{data/fig_fixed_point.dat}\figfixedpoint;
\begin{groupplot}[
	group style = {group size = 2 by 1, horizontal sep = 4em, vertical sep = 5em}, 
	width=0.5\textwidth,
	height = 0.35\textwidth,
	ymode = log, 
	xmin = 0, 
	xmax = 20.5,
	xtick = {1,5,10,15,...,50},
	xlabel = {degree $(n-1,n)$},
	]
	\nextgroupplot[title = Condition Number at Termination,
		ymin = 1,
		ymax = 1e18,
		ytickten = {0,3,6,9,...,18},
		ylabel = {condition number},
		legend columns = -1,
		legend entries = {SK iteration, Vector Fitting},
		legend to name = legend_fixed,
	]
	\addplot[mark=x, only marks, colorbrewerA3] table
		[x=n, y=cond_sk] {\figfixedpoint};
	\addplot[mark=+, only marks, colorbrewerA4] table
		[x=n, y=cond_vf] {\figfixedpoint};
	\coordinate (c1) at (rel axis cs:0,1);

	\nextgroupplot[title=Residual at Termination,
		ymin = 1e-6,
		ymax = 1e3,
		ytickten = {-12,-9,...,-3,0,3},
		ylabel = {$ (\| \ve r\|_2 -\|\ve r_{LS}\|)/\|\ve r_{LS}\|$},
		ylabel shift = -7pt,
		]
	\addplot[mark=x, only marks, colorbrewerA3] table
		[x=n, y expr = \thisrow{res_sk}/\thisrow{res_pr} - 1] {\figfixedpoint};
	\addplot[mark=+, only marks, colorbrewerA4] table
		[x=n, y expr = \thisrow{res_vf}/\thisrow{res_pr} - 1] {\figfixedpoint};
	\coordinate (c2) at (rel axis cs:1,1);
	
    \coordinate (c3) at ($(c1)!.5!(c2)$);
    \node[below] at (c3 |- current bounding box.south)
      {\pgfplotslegendfromname{legend_fixed}};
\end{groupplot}
\end{tikzpicture}
\caption[Conditioning for Fixed Point Methods]{
	Here we use the CD player model described in \cref{fig:local_minima}
	to exhibit the conditioning number of the iteration matrix used by 
	the SK iteration \cref{eq:sk_step} and Vector Fitting \cref{eq:vecfit} with $b_0=1$
	in the left plot.
	Although based on similar principles, the condition number of the Vector Fitting
	step remains well conditioned despite increasing degree.
	The right plot illustrates the 
	difference between the residual norm found via optimization
	and the residual norm found via these two fixed point methods starting at the optimizer.
	Although the difference is slight,
	the nonlinear least squares approach almost always yields a smaller residual norm.
}
\label{fig:fixed_point}
\end{figure}

Although this algorithm frequently converges to fixed points,
the linear system that is solved to update $\ve a$ and $\ve b$ in \cref{eq:skiter}
rapidly becomes ill-conditioned with increasing degree as illustrated in \cref{fig:fixed_point}.
This is not a function of the polynomial basis,
as we have transformed the Legendre basis to be orthogonal on $[-1000i,1000i]$
and the condition number of $\ma \Phi$ and $\ma \Psi$ both remain below 10.
However, the condition number of base iteration matrix
$\begin{bmatrix} -\ma \Phi & \diag(\ve f)[\ma \Psi]_{\cdot,1:n}\end{bmatrix}$
and that of the scaling $\diag(\ma \Psi\ve b)^{-1}$ grow rapidly, 
combining in the large condition number seen.
This motivates adaptive basis used by Vector Fitting.

\subsubsection{Vector Fitting}
Vector Fitting uses a similar approach to the SK iteration,
but makes two subtle changes, which combined remove the ill-conditioning present in the SK iteration.
In the following discussion we will assume we are constructing a rational approximation of degree $(n-1,n)$,
and later show how to increase the numerator degree for any $m>n-1$. 

The first difference with the SK iteration is that vector fitting uses
a Lagrange basis with nodes $\ve \lambda^{(\ell)}$ that change at each iteration:
\begin{align}
	\phi_k^{(\ell)}(z) &= (z - \lambda_k^{(\ell)})^{-1} \prod_{k=1}^n (z - \lambda_k^{(\ell)}) 
		\quad k = 1,\ldots, n \\ 
	\psi_0^{(\ell)}(z) &= \prod_{k=1}^n (z - \lambda_k^{(\ell)})
	\quad \text{and} \quad  
	\psi_k^{(\ell)} = \phi_k^{(\ell)}(z)
		\quad k = 1,\ldots, n 
\end{align}
Then, as in AAA, we take the ratio and canceling the common product, yielding the parameterization
\begin{equation}
	r(z; \ve a, \ve b,\ve\lambda^{(\ell)}) :=
	 \frac{\sum_{k=1}^n a_k \phi_k^{(\ell)}(z)}{\sum_{k=0}^{n} b_k \psi_k^{(\ell)}(z)}
	 = \frac{ \sum_{k=1}^n a_k (z - \lambda_k^{(\ell)})^{-1} }{b_0 + \sum_{k=1}^n b_k (z - \lambda_k^{(\ell)})^{-1}}.
\end{equation}
As with the SK iteration, we multiply through by the denominator to yield a linear optimization problem,
\begin{equation}\label{eq:vecfit}
	\begin{split}
	\ve a^{(\ell+1)}, \ve b^{(\ell+1)} &= \argmin_{\ve a, \ve b}
		\| \diag(\ve f) \ma \Psi^{(\ell)} \ve b - \ma \Phi^{(\ell)}\ve a\|_2 \\
	\text{where} \quad 
	[\ma \Psi^{(\ell)}]_{j,0} &= 1, \quad 
	[\ma \Psi^{(\ell)}]_{j,k} = (z - \lambda_k^{(\ell)})^{-1}, \quad
	[\ma \Phi^{(\ell)}]_{j,k} = (z - \lambda_k^{(\ell)})^{-1}.
	\end{split}
\end{equation}
As with the SK iteration, we have a free scaling between $\ve a$ and $\ve b$.
In the original version of vector fitting~\cite{GS99}, this was dealt with by fixing $b_0 = 1$.
In an updated version Gustavsen recommends adding a constraint that the average real part of
the denominator is one~\cite[eq.~(8)]{Gus06}: 
\begin{equation}
	\sum_{j=1}^N \sum_{k=1}^n \Re[ b_0 + b_k (z_j - \lambda_k^{(\ell)})^{-1} ] = N. 
\end{equation} 
Unlike the SK iteration, vector fitting iteration~\cref{eq:vecfit}
does not include a scaling by the previous denominator.
Instead, after every iteration the Lagrange nodes $\ve \lambda^{(\ell)}$
are updated to be the roots the denominator polynomial by solving the eigenvalue problem~\cite[eq.~(5)]{Gus06}:
\begin{equation}
	\ma A \ve x_k = \lambda_k^{(\ell+1)} \ve x_k, 
	\quad \ma A = \diag(\ve\lambda^{(\ell)}) - \ve 1 \ve b_{1:n}^\trans/b_0 \in \C^{n\times n}.
\end{equation}
Then, if $\ve \lambda^{(\ell)}$ converges to $\ve \lambda^{(*)}$,
the denominator coefficients $\ve b$ converge to $\ve e_0$,
and the denominator polynomial $q$ converges to one,
and the error committed by multiplying by the denominator vanishes.
Moreover, as $q\to 1$, we recover a pole-residue expansion of $r$
with poles $\ve\lambda^{(\ell)}$ and residues $\ve a$.
To extend this algorithm for numerators of degree $m > n - 1$, 
it is sufficient to append columns to $\ma \Phi$.
The complete algorithm is given in \cref{alg:vecfit}.

\begin{algorithm}[t]
\begin{minipage}{\linewidth}
\begin{algorithm2e}[H]
	\Input{Input-output pairs $\lbrace z_j, f(z_j)\rbrace_{j=1}^N$, initial poles $\ve\lambda_0$, degree $(m,n)$}
	\Output{Rational approximation $r(z) = \sum_{k=1}^n \rho_k(z-\lambda)^{-1} + \sum_{k=1}^m c_k \varphi_k(z)$}
	\For{$\ell=0,1,\ldots$}{
		Form Cauchy matrix $[\ma C^{(\ell)}]_{j,k} = (z_j - \lambda_k^{(\ell)})^{-1}$\;
		Form additional polynomial basis: $\ma \Phi\in \C^{N\times (m-n)}$, $[\ma \Phi]_{j,k} = \varphi_k(z_j)$\;
		Solve 
		$\displaystyle
		\begin{bmatrix}
			\ve a^{(\ell+1)} \\
			\ve b^{(\ell+1)}
		\end{bmatrix}
		\leftarrow \argmin_{\ve x}
			\left\|
			\ve f - 
			\begin{bmatrix}
				-\ma C^{(\ell)} & -\ma \Phi &\diag(\ve f) \ma C^{(\ell)} 
			\end{bmatrix}
			\ve x
			\right\|_2
		$\;
		Set $\ve\lambda^{(\ell+1)}$ to be eigenvalues of $\diag(\ve\lambda^{(\ell)}) - \ve 1 \ve b^\trans$\;
		\lIf{$\|\ve b\|_2<$ tol}{{\bf break}}
	}
	Residues $\ve \rho \leftarrow [\ve a^{(*)}]_{1:n}$\;
	Coefficients of polynomial terms $\ve c \leftarrow [\ve a^{(*)}]_{n+1:m+1}$\;
	Poles $\ve\lambda \leftarrow \ve\lambda^{(*)}$\;
\end{algorithm2e}
\vspace{-1em}
\end{minipage}
\caption{Vector Fitting}
\label{alg:vecfit}
\end{algorithm}

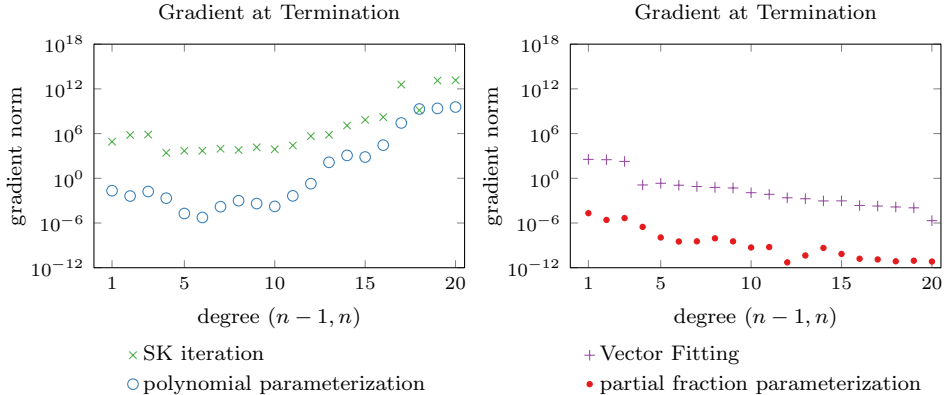
\begin{figure}
\centering
\begin{tikzpicture}
\pgfplotstableread{data/fig_fixed_point.dat}\figfixedpoint;
\begin{groupplot}[
	group style = {group size = 2 by 1, horizontal sep = 4em, vertical sep = 5em}, 
	width=0.5\textwidth,
	height = 0.35\textwidth,
	ymode = log, 
	xmin = 0, 
	xmax = 20.5,
	xtick = {1,5,10,15,...,50},
	xlabel = {degree $(n-1,n)$},
	legend columns = 1,
	legend style = {
		at = {(0.5, -0.3)},
		anchor = north,
	},
	]

	\nextgroupplot[title=Gradient at Termination,
		ymin=1e-12,
		ymax = 1e18,
		ytickten={-12,-6,..., 18},
		ylabel = {gradient norm},
		ylabel shift = -5pt,
		legend entries = {SK iteration, polynomial parameterization},
	]
	\addplot[mark=x, only marks, colorbrewerA3] table
		[x=n, y = grad_sk] {\figfixedpoint};
	\addplot[mark=o, only marks, colorbrewerA2, mark size = 2] table
		[x=n, y = grad_pb] {\figfixedpoint};
	
	\nextgroupplot[title=Gradient at Termination,
		ymin=1e-12,
		ymax = 1e18,
		ytickten={-12,-6,..., 18},
		ylabel = {gradient norm},
		ylabel shift = -5pt,
		legend entries = {Vector Fitting, partial fraction parameterization},
	]
	\addplot[mark=+, only marks, colorbrewerA4] table
		[x=n, y = grad_vf] {\figfixedpoint};
	\addplot[mark=*, only marks, colorbrewerA1, mark size = 1] table
		[x=n, y = grad_pr] {\figfixedpoint};

\end{groupplot}
\end{tikzpicture}
\caption[Gradient at Termination]{
	As an illustration that fixed point methods do not satisfy the first order optimality criteria to high precision
	on the CD player model of \cref{fig:local_minima}.
	These plots show the norm of the gradient at termination after starting from the least squares optimum
	compared with the gradient norm of the comparable nonlinear least squares parameterization.
	Here we see that the optimization based approaches have a gradient norm several orders of magnitude smaller
	than these fixed point methods.
}
\label{fig:grad}
\end{figure}

Although the iterates of this algorithm are substantially better conditioned than
those of the SK iteration, as illustrated in \cref{fig:fixed_point},
this algorithm is not without its concerns.
There are examples of input data where all fixed points are repelling, 
causing the algorithm to iterate indefinitely~\cite{LA13},
which Lefteriu and Antoulas suggest fixing by adding a Newton step.
However, this does not address a more subtle issue:
although the fixed points of this algorithm often provide excellent 
rational approximations, as evidenced by \cref{fig:fixed_point},
these fixed points do not satisfy the first order necessary conditions 
for the least squares rational approximation problem~\cref{eq:ls_ratfit}.
Instead, even when initialized at the least squares optimizer,
the gradient associated with the fixed point of both algorithms 
was substantially larger than that generated using optimization, 
as illustrated in \cref{fig:grad}.
This does not appear to be a numerical artifact
because when initialized at the least squares optimizer,
the first step changed $\ve b$ by at least $10^{-3}$
in both the SK iteration and Vector Fitting for every $n$ test.
This also helps motivate our development of nonlinear least squares approaches in the next two sections.

\section{Optimization Using a Polynomial Basis\label{sec:polyopt}}
An alternative to both Loewner framework and the fixed point iterations 
presented previously is to consider the rational approximation problem
in the light of standard optimization algorithms for least squares problems~\cite{HPS13}.
In this section we will discuss how to apply these results
when the rational approximant is parameterized as a ratio of polynomials
\begin{equation}
	r(z; \ve a, \ve b) := 
		\frac{p(z;\ve a)}{q(z; \ve b)} = \frac{\sum_{k=0}^m a_k \phi_k(z)}{\sum_{k=0}^n b_k \psi_k(z)}
\end{equation}
where $\lbrace \phi_k \rbrace_{k=0}^m$ and $\lbrace \psi_k \rbrace_{k=0}^n$ are two polynomial bases.
This approach has the advantage of being able to represent any degree rational approximant,
whereas the pole-residue parameterization described in the next section
is restricted to degree $(m,n)$ where $m\ge n-1$.
Unfortunately as with the SK-iteration,
the use of a polynomial basis makes this method ill-conditioned 
and of limited utility for rational approximations of moderate dimension, e.g., $n\approx 20$.
However, the ideas developed in this approach are later applied in the next section
to construct a real rational approximation.
In this section, we first derive how to use Variable Projection (VARPRO)~\cite{GP73}
to construct an optimization problem over $\ve b$ alone
and then discuss how to construct a real rational approximation.

\subsection{Variable Projection}
To apply Variable Projection to this rational approximation problem, 
we first state the optimization problem in terms of two Vandermonde matrices
$\ma \Phi\in \C^{N\times (m+1)}$ and $\ma \Psi \in \C^{N\times (n+1)}$:
\begin{equation}
\ma \Phi = \begin{bmatrix} \phi_0(z_1) & \cdots & \phi_m(z_1) \\
				\vdots & & \vdots \\
				\phi_0(z_M) & \cdots & \phi_m(z_M)
			\end{bmatrix},
\quad 
\ma \Psi = \begin{bmatrix} \psi_0(z_1) & \cdots & \psi_m(z_1) \\
				\vdots & & \vdots \\
				\psi_0(z_M) & \cdots & \psi_m(z_M)
			\end{bmatrix},
\end{equation}
as then the rational approximation is
\begin{equation}\label{eq:poly_pre_vp}
	\minimize_{\ve a \in \C^{m+1}, \ve b \in \C^{n+1}}
		\| \ma W[ \ve f - \diag(\ma \Psi \ve b)^{-1} \ma \Phi \ve a ]\|_2.
\end{equation}
The key insight in Variable Projection is that if the nonlinear parameter $\ve b$
is held fixed, $\ve a$ can be written in terms of the pseudoinverse, denoted $^+$:
\begin{equation}
	\ve a = [ \ma W \diag(\ma \Psi \ve b)^{-1}\ma \Phi)]^+ \ma W \ve f.
\end{equation}
If we substitute this value of $\ve a$ into~\cref{eq:poly_pre_vp}, 
we recover an equivalent optimization problem over $\ve b$ alone:
\begin{equation}\label{eq:poly_vp}
	\minimize_{\ve b \in \C^{n+1}}
		\| \ma P_{\ma \Omega(\ve b)}^\perp \ma W \ve f\|_2,
	\quad
	\ma \Omega(\ve b) := \ma W \diag(\ma \Psi \ve b)^{-1} \ma \Phi,
\end{equation}
and where $\ma P_{\ma \Omega(\ve b)}^\perp$ is the projector onto 
the orthogonal complement of the range of $\ma \Omega(\ve b)$.
Defining the interior of this minimization problem as the residual
$\ve r: \C^{n+1} \to \C^N$, 
\begin{equation}
	\ve r(\ve b) := \ma P_{\ma \Omega(\ve b)}^\perp \ma W \ve f,
\end{equation}
Golub and Pereyra then give a formula the Jacobian of $\ve r(\ve b)$
where
\begin{equation}\label{eq:poly_vp_jac}
		\frac{\partial \ve r(\ve b)}{\partial b_k} = 
		- \ma P_{\ma \Omega(\ve b)}^\perp
		\frac{\partial \ma \Omega(\ve b)}{\partial b_k}
		\ma \Omega(\ve b)^+\ma W\ve f
		-
		\ma \Omega(\ve b)^{+*} 
		\frac{\partial \ma \Omega(\ve b)^*}{\partial \conj{b_k} }
		\ma P_{\ma\Omega(\ve b)}^\perp \ma W\ve f
\end{equation}
where we have invoked Wirtinger calculus~\cite[App.~2]{SS10} to extend this result.
However, we must be careful here as $\ve r(\ve b)$ is not an analytic function of $\ve b$.
Instead we will define the Jacobian of $\ve r(\ve b)$ split into real and imaginary parts:
\begin{equation}
	\uve r(\ve b) := \begin{bmatrix} \Re \ve r(\ve b) \\ \Im \ve r(\ve b)\end{bmatrix}
\end{equation}
and construct the Jacobian with respect to $\ve b$ also split into real and imaginary parts:
\begin{equation}
	\uma J(\ve b) := \begin{bmatrix}
				\frac{\partial \Re \ve r(\ve b)}{\partial \Re \ve b} &
				\frac{\partial \Re \ve r(\ve b)}{\partial \Im \ve b} \\ 
				\frac{\partial \Im \ve r(\ve b)}{\partial \Re \ve b} &
				\frac{\partial \Im \ve r(\ve b)}{\partial \Im \ve b}  
			\end{bmatrix}.
\end{equation}l
Then we define matrices $\ma K(\ve b)$ and $\ma L(\ve b)$ 
related the two terms in the Jacobian~\cref{eq:poly_vp_jac}
\begin{equation}\label{eq:poly_KM}
	[\ma K(\ve b)]_{\cdot,k} \!:=\!  
		- \ma P_{\ma \Omega(\ve b)}^\perp
		\frac{\partial \ma \Omega(\ve b)}{\partial \Re b_k}
		\ma \Omega(\ve b)^+\ma W\ve f,
	\ \
	[\ma L(\ve b)]_{\cdot, k} \!:=\!
		-\ma \Omega(\ve b)^{+*} 
		\frac{\partial \ma \Omega(\ve b)^*}{\partial \Re b_k }
		\ma P_{\ma\Omega(\ve b)}^\perp \ma W\ve f, 
\end{equation} 
where the derivative of $\ma \Omega(\ve b)$ is 
\begin{align}
	\frac{\partial \ma \Omega(\ve b) }{\partial \Re b_k} =&
		-\ma W\diag(\ma \Psi \ve b)^{-1} \frac{\partial \diag(\ma \Psi \ve b)}{\partial \Re b_k}
		\diag(\ma \Psi \ve b)^{-1} \ma \Phi \\
		=& -\ma W\diag(\ma \Psi \ve e_k) \diag(\ma \Psi \ve b)^{-2} \ma \Phi.
\end{align}
Using these two matrices we note that by the chain rule, 
derivatives with respect to $\Im b_k$ simply multiply
$\ma K(\ve b)$ and $\ma L(\ve b)$ by $i$.
Hence the Jacobian $\uma J(\ve b)$ is
\begin{equation}\label{eq:poly_JRI}
	\uma J(\ve b) = \begin{bmatrix}
		\Re \ma K(\ve b) + \Re \ma L(\ve b) & -\Im \ma K(\ve b) + \Im \ma L(\ve b)\\
		\Im \ma K(\ve b) + \Im \ma L(\ve b) & \phantom{-}\Re \ma K(\ve b) - \Re \ma L(\ve b)
	\end{bmatrix},
\end{equation}
recalling that the part of the Jacobian corresponding to $\ma L$ 
appears with a derivative with respect to $\conj b_k$.
An algorithm to construct this residual and Jacobian is given in \cref{alg:polyJ}.

\begin{algorithm}[t]
\begin{minipage}{\linewidth}
\begin{algorithm2e}[H]
	\Input{Polynomial coefficients $\ve b\in \C^{n+1}$}
	\Output{Residual $\uve r\in \R^{2N}$ and Jacobian $\uma J\in \R^{(2N)\times (2n + 2)}$}
	Form $\ma \Omega$ from \cref{eq:poly_vp}: $\ma \Omega \leftarrow \ma W \diag(\ma \Psi \ve b)^{-1} \Phi$\;
	Compute short form QR decomposition: $\ma Q \ma R \leftarrow \ma \Omega$\;
	$\ve r \leftarrow \ma W \ve f - \ma Q \ma Q^* \ma W \ve f$\;
	$\uve r \leftarrow \begin{bmatrix} \Re \ve r \\ \Im \ve r\end{bmatrix}$\;
	$\ve a \leftarrow  \ma R^+\ma Q^* \ma W \ve f$\; 
	\label{alg:polyJ:a}
	$[\ma K]_{\cdot,k} \leftarrow  [\ma I - \ma Q \ma Q^*]\ma W\diag(\ma \Psi \ve e_k)\diag(\ma \Psi \ve b)^{-2} \ma \Phi \ve a
	\quad k=0,\ldots,n$\;
	$[\ma L]_{\cdot,k} \leftarrow \ma Q\ma R^{+*} \ma \Phi^*\diag(\ma \Psi \ve e_k)^*\diag(\ma \Psi \ve b)^{-2*} \ma W^*\ve r
	\quad k =0,\ldots, n$\;
	$\uma J \leftarrow \begin{bmatrix}
			\Re \ma K + \Re \ma L & -\Im \ma K + \Im \ma L\\
			\Im \ma K + \Im \ma L & \phantom{-}\Re \ma K - \Re \ma L
		\end{bmatrix}
	$\;  \label{alg:polyJ:J}
\end{algorithm2e}
\vspace{-1em}
\end{minipage}
\caption{Residual and Jacobian for Complex Polynomial Parameterization}
\label{alg:polyJ}
\end{algorithm}

\subsection{Real Rational Approximation\label{sec:polyopt:real}}
If we wish to impose the constraint the approximant $r$ is a real rational function, $r \in \R_{m,n}(\R)$,
one approach is modify the parameterization such that $r$ is necessarily in this class.
Using the polynomial basis, 
we can enforce this constraint by choosing the bases $\lbrace \phi_k\rbrace_{k=0}^m$ and $\lbrace \psi_k\rbrace_{k=0}^n$
consist solely of real polynomials (polynomials with only real coefficients),
and then requiring the coefficients $\ve a$ and $\ve b$ to both be real.
This causes several changes to our formula for the Jacobian.
First, as $\ve b$ is real, the Jacobian now only has $n+1$ columns.
Next, we need to ensure that the solution for $\ve a$ is real
as well as when the projector $\ma P_{\ma \Omega(\ve b)}$ is used.
To do so, we instead use the projector from $\ma \Omega(\ve b)$ split into real and imaginary parts $\ma P_{\uma \Omega(\ve b)}$:
\begin{equation}
	\uma \Omega(\ve b) := \begin{bmatrix} \Re \ma \Omega(\ve b) \\ \Im \ma \Omega(\ve b) \end{bmatrix}.
\end{equation}
Then we form the Jacobian using this projector as before,
but instead using $\ma W\ve f$ split into real and imaginary parts:
\begin{equation}
	[\uma J(\ve b)]_{\cdot,k} = 
	-\ma P_{\uma \Omega(\ve b)}^\perp \frac{\partial \uma \Omega(\ve b)}{\partial \Re b_k} \uma \Omega(\ve b)^+\uve f
		- \uma \Omega(\ve b)^{+*} \frac{\partial \uma \Omega^*}{\partial \Re b_k} \ma P_{\uma \Omega(\ve b)}^\perp \uve f;
	\quad \uve f = \begin{bmatrix} \Re \ma W \ve f \\ \Im \ma W \ve f\end{bmatrix}. 
\end{equation}
These modifications are summarized in \cref{alg:polyJ_real}.
\begin{algorithm}[t]
\begin{minipage}{\linewidth}
\begin{algorithm2e}[H]
	\Input{Polynomial coefficients $\ve b\in \R^{n+1}$}
	\Output{Residual $\uve r\in \R^{2N}$ and Jacobian $\uma J\in \R^{(2N)\times (n + 1)}$}
	Form $\ma \Omega$ from \cref{eq:poly_vp}: $\ma \Omega \leftarrow \ma W \diag(\ma \Psi \ve b)^{-1} \Phi$\;
	$\uma \Omega \leftarrow \begin{bmatrix} \Re \ma \Omega \\ \Im \ma \Omega \end{bmatrix}\qquad
	\uma f \leftarrow \begin{bmatrix} \Re \ma W\ve f \\ \Im \ma W \ve f \end{bmatrix}$
	\;
	Compute short form QR decomposition: $\uma Q \uma R \leftarrow \uma \Omega$\;
	$\uve r \leftarrow \uma f - \uma Q\, \uma Q^\trans \uma f\qquad \ve r \leftarrow [\uve r]_{1:N} + i[\uve r]_{N+1:2N}$\;
	$\ve a \leftarrow  \uma R^+\uma Q^\trans \uma f$\;
		$[\uma K]_{\cdot,k} 
			\leftarrow  [\ma I - \uma Q\, \uma Q^\trans] 
			\begin{bmatrix}
			 \Re \ma W\diag(\ma \Psi \ve e_k)\diag(\ma \Psi \ve b)^{-2} \ma \Phi \ve a \\
			 \Im \ma W\diag(\ma \Psi \ve e_k)\diag(\ma \Psi \ve b)^{-2} \ma \Phi \ve a 
			\end{bmatrix}
			\quad k=0,\ldots,n$\;
		$[\uma L]_{\cdot,k} 
			\leftarrow \uma Q\uma R^{+\trans} 
			\begin{bmatrix}
				 \Re \ma \Phi^*\diag(\ma \Psi \ve e_k)^*\diag(\ma \Psi \ve b)^{-2*} \ma W^*\ve r \\
				 \Im \ma \Phi^*\diag(\ma \Psi \ve e_k)^*\diag(\ma \Psi \ve b)^{-2*}\ma W^*\ve r
			\end{bmatrix}
			\quad k=0,\ldots,n$\;
	$\uma J \leftarrow \uma K + \uma L$\; 
\end{algorithm2e}
\vspace{-1em}
\end{minipage}
\caption{Residual and Jacobian for Real Polynomial Parameterization}
\label{alg:polyJ_real}
\end{algorithm}

\subsection{Normalization and Conditioning}
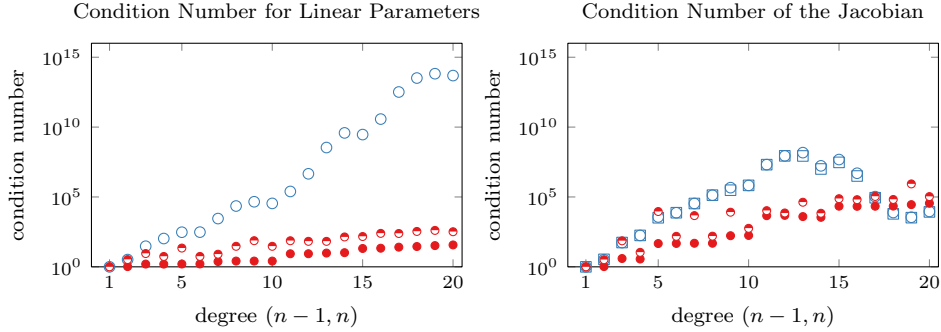
\begin{figure}

\begin{tikzpicture}
\pgfplotstableread{data/fig_condition.dat}\figcondition;
\begin{groupplot}[
	group style = {group size = 2 by 1, horizontal sep = 4em}, 
	width=0.5\textwidth,
	height = 0.35\textwidth,
	ymode = log, 
	xmin = 0, 
	xmax = 20.5,
	xtick = {1,5,10,15,...,50},
	xlabel = {degree $(n-1,n)$},
	ymin = 1,
	ymax = 1e16,
	ytickten = {0,5,10,15},
	ylabel = condition number
	]
	\nextgroupplot[title=Condition Number for Linear Parameters]
	\addplot[mark=o, only marks, colorbrewerA2] table
		[x=n, y=cond_omega] {\figcondition};
	\addplot[mark=*, mark size = 1.5, only marks, colorbrewerA1] table
		[x=n, y=cond_V] {\figcondition};
	
	\addplot[mark=halfcircle*, mark size = 1.5, only marks, colorbrewerA1] table
		[x=n, y=cond_V_real] {\figcondition};
	
	\nextgroupplot[title=Condition Number of the Jacobian]
	\addplot[mark=o, only marks, colorbrewerA2] table
		[x=n, y=cond_poly_vp] {\figcondition};
	\addplot[mark=square, only marks, colorbrewerA2] table
		[x=n, y=cond_poly_vp_grass] {\figcondition};
	\addplot[mark=*, mark size = 1.5, only marks, colorbrewerA1] table
		[x=n, y=cond_pr_vp] {\figcondition};
	\addplot[mark=halfcircle*, mark size = 1.5, only marks, colorbrewerA1] table
		[x=n, y=cond_pr_vp_real] {\figcondition};
\end{groupplot}
\end{tikzpicture}
\caption[Condition numbers of basis]{
	Here we show two sets of condition numbers relevant to optimization
	in both the polynomial and pole-residue bases.
	The left plot shows the condition number associated with finding the linear parameters.
	For the polynomial basis, this is the condition number of $\ma \Omega(\ve b)$ given in \cref{eq:poly_vp}
	denoted by \tikz{\draw[colorbrewerA2] plot[mark=o] (0,0);};
	for the pole-residue basis, this is the condition number of $\ma \Lambda(\ve\lambda)$ given in \cref{eq:pr_Lam}
	denoted by \tikz{\draw[colorbrewerA1] plot[mark=*, mark size = 1.5] (0,0);};
	for the quadratic partial fraction expansion used in \cref{sec:propt:real}
	this is $\ma \Theta(\ve b)$ given by~\cref{eq:Theta}
	denoted by \tikz{\draw[colorbrewerA1] plot[mark=halfcircle*, mark size = 1.5] (0,0);};
	The right plot shows the condition number of the Jacobian
	with the polynomial basis using a monic constraint shown as \tikz{\draw[colorbrewerA2] plot[mark=o] (0,0);},
	and norm constraint shown as \tikz{\draw[colorbrewerA2] plot[mark=square] (0,0);};
	the condition number of the pole-residue Jacobian is shown as 
	\tikz{\draw[colorbrewerA1] plot[mark=*, mark size = 1.5] (0,0);} for the complex case
	and \tikz{\draw[colorbrewerA1] plot[mark=halfcircle*, mark size = 1.5] (0,0);} for the real case.
	In each case these matrices are evaluated at a local optimum of the CD player model described in \cref{fig:local_minima}.	
}
\label{fig:condition}
\end{figure}

As with the SK-iteration, 
we now face a choice of how to remove the additional degree of freedom in our choice of $\ve b$.
One simple option would be to simply fix one of the entries;
for example, forcing $b_n=1$ which then yields small Jacobian with only $2n$ columns.
Another option is to constrain the norm of $\ve b$, e.g., $\| \ve b\|_2 = 1$,
which is more numerically sound if the coefficient we fixed is small or zero at the optimum.
However, fixing the norm of $\ve b$ is substantially more difficult to implement;
for example one approach would be perform optimization on the Grassmann manifold~\cite{EAS98}.
Unfortunately, neither approach is ultimately helpful.
As shown in \cref{fig:condition}, 
the condition number of $\ma \Omega(\ve b)$ grows rapidly with increasing degree.
Similarly, the Jacobian is increasingly ill-conditioned until the loss of precision in the pseudoinverse of 
$\ma \Omega(\ve b)$ causes the condition number to artificially decrease;
using the full Jacobian, ignoring its two dimensional nullspace, does not fix the conditioning issues either.
This motivates using the pole-residue basis we discuss in the next section.

\section{Optimization Using a Partial Fraction Parameterization\label{sec:propt}}
An alternative to the ratio of two polynomials used in the previous section
is to instead consider the partial fraction expansion of $r$ into a sum of
degree $(0,1)$ rational functions.
For a degree $(m,n)$ rational function where $m \ge n-1$
and $q$ has $n$ distinct roots $\lbrace \lambda_k \rbrace_{k=1}^n$ we can write
\begin{equation}
	r(z; \ve a, \ve b) := 
		\frac{p(z; \ve a)}{q(z;\ve b)} = 
		\frac{\sum_{k=0}^m a_k \phi_k(z)}{\sum_{k=0}^n b_k\psi_k(z)} 
	= \sum_{k=1}^n \frac{ \rho_k}{z - \lambda_k} + \sum_{k=0}^{m-n} c_k \varphi_k(z)
	=: r(z; \ve \lambda, \ve \rho, \ve c)
\end{equation}
where $\lbrace\varphi_k\rbrace_{k=0}^{m-n}$ is a basis for polynomials of degree $m-n$.
This partial fraction expansion, also known as a pole-residue expansion,
is much better conditioned than the polynomial basis for rational approximation considered in the previous section,
as evidenced by \cref{fig:condition},
leading to better optimizers of the least squares rational approximation problem.
There is a small price we have paid for this:
the pole-residue basis cannot express higher order poles $(z-\lambda)^2$, $(z-\lambda)^3$, etc.
However, for any polynomial $q$ with multiple roots there is another polynomial $\widetilde q$
arbitrarily close with distinct roots.
Hence for the purposes of rational approximation,
even if $f$ has higher order poles, we can approximate it arbitrarily well in the pole-residue parameterization.

In the following subsections, 
we first derive the VARPRO residual and Jacobian for this problem
and then show how enforce that $r$ is a real rational function.
Here we are forced to modify our parameterization,
choosing a partial fraction expansion into a sum of degree $(1,2)$ rational functions
and we reuse portions of the derivation from~\cref{sec:polyopt:real}.

\subsection{Variable Projection}
To begin, we first write the least squares rational approximation problem using this parameterization
in terms of a Cauchy matrix $\ma C(\ve\lambda) \in \C^{N\times n}$ and a Vandermonde matrix $\ma \Phi\in \C^{N\times (m-n)}$:
\begin{equation}
	\minimize_{\ve \lambda \in \C^n, \ve \rho \in \C^n,\ve c \in \C^{m-n}}
		\| \ma W[ \ve f - \ma C(\ve\lambda) \ve \rho - \ma \Phi \ve c ]\|_2.
\end{equation} 
where
\begin{equation}\label{eq:propt_C_Phi}
\ma C(\ve\lambda)  \!:=\! \begin{bmatrix} (z_1 \!-\! \lambda_1)^{-1}& \cdots & (z_1 \!-\! \lambda_n)^{-1} \\
				\vdots & & \vdots \\
				(z_N \!-\! \lambda_1)^{-1} & \cdots & (z_N \!-\! \lambda_n)^{-1}
			\end{bmatrix},
\ \
\ma \Phi \!:=\! \begin{bmatrix} \varphi_0(z_1) & \cdots & \varphi_m(z_1) \\
				\vdots & & \vdots \\
				\varphi_0(z_M) & \cdots & \varphi_m(z_M)
			\end{bmatrix}.
\end{equation}
As with optimization in the polynomial basis, 
it is helpful to define a single matrix function $\ma \Lambda(\ve \lambda)$,
analogous to $\ma \Omega(\ve b)$ in \cref{eq:poly_vp}: 
\begin{equation}\label{eq:pr_Lam}
	\ma \Lambda( \ve\lambda) := \begin{bmatrix} \ma W\ma C(\ve\lambda) & \ma W\ma \Phi \end{bmatrix},
\end{equation}
leaving the optimization problem:
\begin{equation}
	\minimize_{\ve \lambda \in \C^n, \ve d \in \C^{m}}
		\| \ma W \ve f - \ma \Lambda(\ve\lambda) \ve d\|_2,
		\quad \ve d = \begin{bmatrix} \ve \rho \\ \ve c \end{bmatrix}.
\end{equation}
Then, as before, 
we use VARPRO to convert this into an optimization problem over $\ve \lambda$ alone:
\begin{equation}
	\minimize_{\ve\lambda \in \C^n}
		\| \ma P_{\ma \Lambda(\ve\lambda)}^\perp \ma W \ve f \|_2
	\quad \text{with residual}
	\quad \ve r(\ve\lambda) := \ma P_{\ma \Lambda(\ve \lambda)}^\perp \ma W\ve f.
\end{equation}
Then, defining the two terms in the VARPRO Jacobian $\ma K$ and $\ma L$ analogously to~\cref{eq:poly_KM},
\begin{equation}\label{eq:pr_KM}
	[\ma K(\ve \lambda )]_{\cdot,k} \!:=\!  
		- \ma P_{\ma \Lambda(\ve \lambda)}^\perp
		\frac{\partial \ma \Lambda(\ve \lambda)}{\partial \Re \lambda_k}
		\ma \Lambda(\ve \lambda)^+\ma W\ve f,
	\ \
	[\ma L(\ve \lambda)]_{\cdot, k} \!:=\!
		-\ma \Lambda(\ve \lambda)^{+*} 
		\frac{\partial \ma \Lambda(\ve \lambda)^*}{\partial \Re\lambda_k }
		\ma P_{\ma\Lambda(\ve \lambda)}^\perp \ma W\ve f, 
\end{equation} 
where the derivative $\ma \Lambda(\ve \lambda)$ with respect to $\lambda_k$ is
\begin{equation}
	\frac{\partial \ma \Lambda(\ve \lambda)}{\partial \Re \lambda_k}
		= -\begin{bmatrix} (z_1 - \lambda_k)^{-2} \\ \vdots \\ (z_N - \lambda_k)^{-2} \end{bmatrix}
		\ve e_k^\trans,
\end{equation} 
where $\ve e_k$ is the $k$th unit vector.
Then, using this $\ma K(\ve\lambda)$ and $\ma L(\ve\lambda)$ we can build the 
the Jacobian using~\cref{eq:poly_JRI} as before.
In building the residual and Jacobian in \cref{alg:prJ}, 
we also exploit the rank-1 structure of $\partial \ma \Lambda(\ve\lambda)/\partial \Re\lambda_k$,
to build $\ma K$ and $\ma L$ using matrix-matrix products rather than a loop over $k$
as in \cref{alg:polyJ}.

\begin{algorithm}[t]
\begin{minipage}{\linewidth}
\begin{algorithm2e}[H]
	\Input{Poles $\ve\lambda \in \C^n$}
	\Output{Residual $\uve r\in \R^{2N}$ and Jacobian $\uma J\in R^{(2N)\times (2n)}$}
	Form $ \ma \Lambda \leftarrow \ma \Lambda(\ve\lambda)$ using \cref{eq:pr_Lam}\;
	Compute the short form QR decomposition $\ma Q \ma R \leftarrow \ma \Lambda$\;
	Compute residual $\ve r \leftarrow \ma W \ve f - \ma Q \ma Q^* \ma W \ve f$\;
	$\uve r \leftarrow \begin{bmatrix} \Re \ve r \\ \Im \ve r\end{bmatrix}$\;
	$\begin{bmatrix} \ve \rho \\ \ve c \end{bmatrix}
		\leftarrow \ve d \leftarrow  \ma R^+\ma Q^* \ma W \ve f$\;
	Build $[\ma \Lambda']_{j,k} = -(z_j - \lambda_k)^{-2}$\;
	$\ma K \leftarrow -[\ma I - \ma Q \ma Q^*]\ma W \ma \Lambda' \diag(\ve \rho)$\;
	$\ma L \leftarrow -\ma Q\ma R^{+*}\diag(\ma \Lambda'^*\ve r)$\;
	$\uma J \leftarrow \begin{bmatrix}
			\Re \ma K + \Re \ma L & -\Im \ma K + \Im \ma L\\
			\Im \ma K + \Im \ma L & \phantom{-}\Re \ma K - \Re \ma L
		\end{bmatrix}
	$\;  
\end{algorithm2e}
\vspace{-1em}
\end{minipage}
\caption{Residual and Jacobian for Partial Fraction Parameterization}
\label{alg:prJ}
\end{algorithm}

\subsection{Real Rational Approximation\label{sec:propt:real}}
It is not simple to construct a real rational approximation in a pole-residue basis.
As with the polynomial basis, we can require the basis for the polynomial component
$\lbrace \varphi_k \rbrace_{k=0}^{m-n}$ consist of real polynomials 
and that the coefficients $\ve c$ be real.
However, more complicated constraints are required for the pole-residue porition
to ensure that $r(\conj{z}) = \conj{r(z)}$.
Namely, for every $\lambda_k$ with nonzero imaginary part
that there is another pole $\lambda_{\set I_k}$ that is its complex conjugate pair,
$\lambda_k = \conj{\lambda_{\set I_k}}$
and that this same relationship applies to the residues, $\rho_k = \conj{\rho_{\set I_k}}$. 
Naively enforcing these constraints is not simple:
these pairings can appear if poles leave the real line
and disappear if poles enter the real line.
Alternatively, we could automatically include the conjugate of every non-real $\lambda_k$
but the degree of the rational approximation will change whenever a $\lambda_k$ becomes real.
If we instead try to work with an implicit parameterization of the poles,
such as parametrizing a pair of poles as roots of quadratic polynomial $q(z) =z^2 +b_1 z + b_0$
this parameterization is not differentiable when $q$ has a multiple root.
Instead, here we develop an approach based on a partial fraction expansion of $r$
into a sum of $(1,2)$ rational functions plus the remaining $m-n$ polynomial terms:
\begin{equation}
	r(z; \ve a, \ve b, \ve c ) = 
	\sum_{k=0}^{m-n} c_k \varphi_k(z)
	+
	\begin{cases}
		\displaystyle \phantom{\frac{a_n}{z + b_n}}\phantom{+}
		\sum_{k=1}^{\lfloor n/2 \rfloor} \frac{ a_{2k} z + a_{2k-1}}{z^2 + b_{2k} z + b_{2k-1}}, & n \text{ even;}\\
		\displaystyle \frac{a_n}{z + b_n}
		+ \! \sum_{k=1}^{\lfloor n/2 \rfloor} \frac{ a_{2k} z + a_{2k-1}}{z^2 + b_{2k} z + b_{2k-1}}, 
		& n \text{ odd;}
	\end{cases}
\end{equation}
where $\ve a, \ve b\in \R^{n}$ and $\ve c \in \R^{m-n}$.
This is a middle ground between the pole residue approach we used before and the polynomial basis
considered in \cref{sec:polyopt} that is better conditioned yet still allows us to easily enforce that $r$
is a real rational function by requiring $\ve a$, $\ve b$, and $\ve c$ to be real.

The derivation of the residual and Jacobian for this parameterization largely follows that of \cref{sec:polyopt:real}.  
Defining $\ma \Theta(\ve b)$ in an analogous role to $\ma \Omega(\ve b)$ and $\ma \Lambda(\ve\lambda)$,
we state the optimization problem as
\begin{align}
	&\minimize_{\ve a, \ve b\in \R^n, \ve c \in \R^{m-n}}
		\left\| \begin{bmatrix} \Re \ma W\ve f \\ \Im \ma W \ve f\end{bmatrix}
		 - \begin{bmatrix} \Re \ma \Theta(\ve b) \\ \Im \ma \Theta(\ve b) \end{bmatrix}
		\begin{bmatrix} \ve a \\ \ve c\end{bmatrix} \right\|_2,
		\quad \text{where} \\
	&\ma \Theta(\ve b) := 
	\begin{cases}
		\ma W\begin{bmatrix} \hma \Omega_1([\ve b]_{1:2}) & \cdots & \hma \Omega([\ve b]_{n-1:n}) & \ma \Phi \end{bmatrix}, 
		& n \text{ even};\\
		\ma W\begin{bmatrix} \hma \Omega_1([\ve b]_{1:2}) & \cdots & \hma \Omega([\ve b]_{n-2:n-1}) &
			(\ve z + b_n)^{-1} & \ma \Phi \end{bmatrix}, 
		& n \text{ odd};
	\end{cases}
	\label{eq:Theta}
\end{align}
and where $\hma \Omega(\ve b)$ defined analogously to $\ma \Omega(\ve b)$ 
for a degree $(1,2)$ rational function in the monomial basis with a monic constraint:
\begin{equation}
	\hma \Omega(\ve b) :=
		\diag\left(
			\begin{bmatrix}
				z_1^2 + b_2 z_1 + b_1) \\
				\vdots \\
				z_N^2 + b_2 z_N + b_1 
			\end{bmatrix}
		\right)^{-1} 
		\begin{bmatrix}
			z_1 & 1 \\
			\vdots & \vdots \\
			z_N & 1
		\end{bmatrix}
\end{equation}
and $\ma \Phi$ is defined as in~\cref{eq:propt_C_Phi}.
Then, after applying VARPRO, our residual is
\begin{equation}
	\uve r(\ve b) = \ma P_{\uma \Theta(\ve b)}^\perp \uma f,
	\quad \uma \Theta(\ve b) = \begin{bmatrix} \Re \ma \Theta(\ve b) \\ \Im \ma \Theta (\ve b)\end{bmatrix}
	\quad \uve f = \begin{bmatrix} \Re \ma W \ve f  \\ \Im \ma W \ve f \end{bmatrix}.
\end{equation} 
Then, computing the Jacobian is similar to \cref{sec:polyopt:real},
except now the Jacobian is split into two column blocks, with an additional one column block if $n$ is odd.
This formula is given in \cref{alg:prJ_real}.
Then once we have found an optimum $\ve b$, we can easily convert back to a pole-residue form,
using the quadratic formula to compute the roots of each term and then compute their corresponding residues.

\begin{algorithm}[t]
\begin{minipage}{\linewidth}
\begin{algorithm2e}[H]
	\Input{Parameters for quadratic expansion $\ve b \in  \R^n$}
	\Output{Residual $\uve r\in \R^{2N}$ and Jacobian $\uma J\in \R^{(2N)\times n}$}
	Form $\ma \Theta \leftarrow \ma \Theta(\ve b)$ from \cref{eq:Theta} and 
	$\uma  \Theta \leftarrow \begin{bmatrix} \Re \ma \Theta \\ \Im \ma \Theta \end{bmatrix}$ \;
	Compute the short form QR decomposition $\uma Q \uma R \leftarrow \uma \Theta$\;
	Define $\uve f \leftarrow \begin{bmatrix} \Re \ma W\ve f \\ \Im \ma W \ve f \end{bmatrix}$ \;
	Compute residual $\uve r \leftarrow \uve f - \uma Q \uma Q^* \uve f \quad \ve r \leftarrow [\uve r]_{1:N} + i[\uve r]_{N+1:2N} $\;
	Compute linear coefficients $\begin{bmatrix} \ve a  \\ \ve c \end{bmatrix}
		\leftarrow \ve d \leftarrow  \uma R^+\uma Q^* \uve f$\;
	\For{$k=1,\ldots, \lfloor n/2\rfloor$}{
		$[\uma K]_{\cdot,2k-1} \leftarrow 
		 [\ma I - \uma Q\, \uma Q^*] 
		\begin{bmatrix}
		 \Re \ma W \diag(\ve z^2+b_{2k}\ve z+ b_{2k-1})^{-2} (a_{2k} \ve z + a_{2k-1}) \\
		 \Im \ma W \diag(\ve z^2+b_{2k}\ve z+ b_{2k-1})^{-2} (a_{2k} \ve z + a_{2k-1}) 
		\end{bmatrix}$\;
		$[\uma K]_{\cdot,2k} \leftarrow 
		 [\ma I - \uma Q\, \uma Q^*] 
		\begin{bmatrix}
		 \Re \ma W \diag(\ve z)\diag(\ve z^2+b_{2k}\ve z+ b_{2k-1})^{-2} (a_{2k} \ve z + a_{2k-1}) \\
		 \Im \ma W \diag(\ve z)\diag(\ve z^2+b_{2k}\ve z+ b_{2k-1})^{-2} (a_{2k} \ve z + a_{2k-1}) 
		\end{bmatrix}$\;
		$[\uma L]_{\cdot,2k-1} 
			\leftarrow \uma Q\uma R^{+\trans} 
			\begin{bmatrix}
				 \Re \ma \Phi^*\diag(\ve z^2 + b_{2k} \ve z + b_{2k-1})^{-2*}\ma W^*\ve r \\
				 \Im \ma \Phi^*\diag(\ve z^2 + b_{2k} \ve z + b_{2k-1})^{-2*}\ma W^*\ve r
			\end{bmatrix}$\;
		$[\uma L]_{\cdot,2k} 
			\leftarrow \uma Q\uma R^{+\trans} 
			\begin{bmatrix}
				 \Re \ma \Phi^*\diag(\ve z)\diag(\ve z^2 + b_{2k} \ve z + b_{2k-1})^{-2*}\ma W^*\ve r \\
				 \Im \ma \Phi^*\diag(\ve z)\diag(\ve z^2 + b_{2k} \ve z + b_{2k-1})^{-2*}\ma W^*\ve r
			\end{bmatrix}$\;
	}
	\If{$n$ is odd}{
		$[\uma K]_{\cdot,n} \leftarrow
			[\ma I - \uma Q\, \uma Q^*] 
			\begin{bmatrix}
			 \Re \ma W \diag(\ve z+ b_n)^{-2} a_n \\
		 	\Im \ma W \diag(\ve z + b_n)^{-2} a_n 
		\end{bmatrix}$\;
		$[\uma L]_{\cdot,n} 
			\leftarrow \uma Q\uma R^{+\trans} 
			\begin{bmatrix}
				 \Re \ma \Phi^*\diag(\ve z + b_n)^{-2*}\ma W^*\ve r \\
				 \Im \ma \Phi^*\diag(\ve z + b_n)^{-2*}\ma W^*\ve r
			\end{bmatrix}$\;
	}
	$\uma J \leftarrow \uma K + \uma L$\;  
\end{algorithm2e}
\vspace{-1em}
\end{minipage}
\caption{Residual and Jacobian for Real Partial Fraction Parameterization}
\label{alg:prJ_real}
\end{algorithm}

Although this approach has used elements of the polynomial basis to enforce the real constraint,
the optimization problem has not become substantially more ill-conditioned than the pole-residue approach of the previous section.
As illustrated in \cref{fig:condition}, 
the condition number of this approach stays close to that of the pole-residue basis.
We suspect that this is due to the fact that this case has used a partial fraction expansion
into degree $(1,2)$ rational functions
whereas the pole-residue approach can be interpreted as an expansion into $(0,1)$ rational functions.

\section{Numerical Examples\label{sec:examples}}
As the previous sections have included several examples illustrating the performance 
of our optimization approach using both polynomial and partial fraction bases
for both real and complex rational approximation,
here we focus our attention to features not addressed in previous examples.
In the first example we show the utility of enforcing the real constraint 
and in the second we consider a rational approximation problem
with a nontrivial weight matrix $\ma W$ related to the projected $\Htwo$ model reduction problem.

\subsection{Employing the Real Constraint}
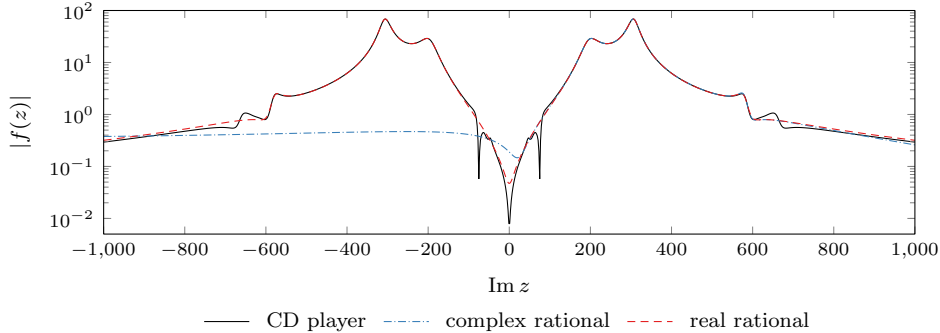
\begin{figure}
\begin{tikzpicture}
\centering
\begin{axis}[
	ymode = log,
	width = 0.95 \textwidth,
	height = 0.35 \textwidth,
	unbounded coords=discard,
	ylabel = $|f(z)|$,
	xlabel = $\Im z$,
	clip mode=individual,
	ymin = 5e-3,
	ymax = 1e2,
	ytickten = {-2,-1,0,1,2},
	legend entries = {CD player, complex rational, real rational },
	legend columns = -1,
	legend style = { 
		at = {(0.5, -0.3)}, 
		anchor = north,	
		column sep=1ex,
	}	
	]
	\addplot[black] plot[x=z_imag, y = h] table {data/fig_real_h.dat};
	\addplot[colorbrewerA2, densely dashdotted] plot[x=z_imag, y = h_complex] table {data/fig_real_h_complex.dat};
	\addplot[colorbrewerA1, densely dashed] plot [x=z_imag, y = h_real] table {data/fig_real_h_real.dat};
\end{axis}
\end{tikzpicture}
\caption{The samples of the CD player model described in \cref{fig:local_minima}
	along with degree $(5,6)$ rational fits using a partial fraction expansion with and without the real constraint
	when fitting to samples with $\Im z >0$.
}
\label{fig:real}
\end{figure}

As an example of how requiring the real constraint can prove beneficial,
\cref{fig:real} illustrates the rational fits to only samples with positive imaginary part.
Without imposing the real constraint, the resulting approximation only fits 
where there are samples with positive part.
However by adding the constraint that $r$ be real, 
the resulting approximation does equally well on both halves of the data.

\subsection{Application to Projected $\Htwo$ Model Reduction\label{sec:examples:ph2}}
In recent work by the authors, 
the $\Htwo$ norm is approximated by its projection onto a finite dimensional subspace~\cite{HM18x}.
Then, as a step in the model reduction process, 
it is necessary to construct a real least squares rational approximation of degree $(n-1,n)$ in a weighted norm:
\begin{equation}
	\minimize_{r \in \set R_{n-1,n}(\R)} \| \ma W[f(\set Z) - r(\set Z)]\|_2 
\end{equation}
where the weight matrix $\ma W$ is the inverse square root of a Cauchy matrix $\ma M(\ve z)$:
\begin{equation}
	\ma W = \ma M(\ve z)^{-1/2}, \qquad 
		\ma M(\ve z) = 
		\begin{bmatrix} 
			(z_1 + \conj{z_1})^{-1} & \cdots & (z_1 + \conj{z_N})^{-1} \\
			\vdots & & \vdots \\
			(z_N + \conj{z_1})^{-1} & \cdots & (z_N + \conj{z_N})^{-1}
		\end{bmatrix}
\end{equation}
and each $z_j$ is in the right half plane.
Using Demmel's results~\cite{Dem00}, 
we can compute the SVD of $\ma M(\ve z) = \ma U \ma \Sigma \ma V^*$ to high relative accuracy
and set $\ma W = \ma \Sigma^{-1/2}\ma U^*$.

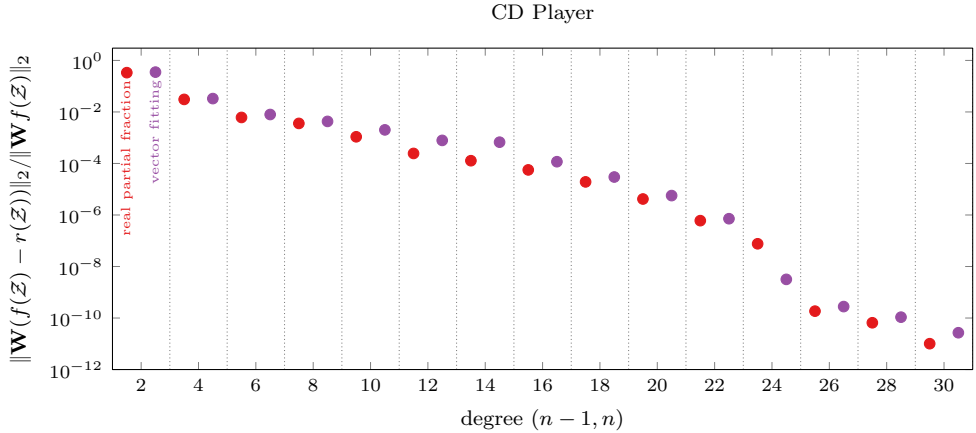
\begin{figure}

\begin{tikzpicture}
\begin{axis}[
	width = \textwidth,
	height = 0.45\textwidth,
	ymode = log,
	xmin = 1,
	xmax = 31,
	ymax = 3e0,
	ymin = 1e-12,
	ytickten = {-14,-12,...,0},
	xlabel = {degree $(n-1,n)$},
	ylabel = {$\| \ma W(f(\set Z) - r(\set Z))\|_2/\|\ma W f(\set Z)\|_2$},
	title = {CD Player},	
	]

	\addplot[colorbrewerA1, mark = *, only marks, mark size = 2, ] table
		[x expr = \thisrow{n} -0.5, y = res_pr]
		{data/fig_ph2.dat};
	
	\addplot[colorbrewerA4, mark = *, only marks, mark size = 2, ] table
		[x expr = \thisrow{n} + 0.5, y = res]
		{data/fig_ph2_matlab.dat};
	
	\foreach \rn in {1,3,...,31}{
		\addplot[gray, densely dotted] coordinates {(\rn,1e-12) (\rn, 3)};	
	}
	
	\draw [colorbrewerA1] (1.5, 0.4) node [anchor = east, rotate = 90] {\tiny real partial fraction};
	\draw [colorbrewerA4] (2.5, 0.4) node [anchor = east, rotate = 90] {\tiny vector fitting};

\end{axis}
\end{tikzpicture}
\caption[Null]{
	The residual norm on the weighted rational fitting problem emerging from the projected $\Htwo$ approach
	as described in \cref{sec:examples:ph2}
	using the real partial fraction expansion described in \cref{sec:propt:real}
	and the {\tt vectfit3} implementation of vector fitting which does not support non-diagonal weightings.
}
\label{fig:ph2}
\end{figure}

As an example of this weighted problem, \cref{fig:ph2} considers the CD player model
sampled a sequence of points with imaginary part uniformly sampled between $-1000i$ and $1000i$:
80 points with $\Re z=0.001$, 40 points with $\Re z=0.01$, 20 points with $\Re z=0.1$,
and 10 points with $\Re z = 1$.
As expected, the {\tt vectfit3} implementation of vector fitting 
does worse in almost every case as it does not support non-diagonal weight matrices.

\section{Discussion\label{sec:discussion}}
Here we have shown how to construct least squares rational approximants 
using standard optimization techniques in two bases, 
a ratio of polynomials and a more stable partial fraction expansion,
and with and without constraining the rational approximant to be real.
Although this optimization approach often finds spurious local minima, 
we have found that the AAA algorithm provides an effective initialization.
Moreover, unlike existing approaches, 
we are able to find rational approximants that statisfy the first order necessary conditions
to high precision.

Our focus has been on the scalar rational approximation problem,
however many applications in systems theory
seek a rational function with vector or matrix valued output:
\begin{equation}
	\set R_{m,n}(\C, \C^{s\times t}) := 
		\lbrace q(z)^{-1} \ma P(z) | q \in \set P_n(\C), \ \ma P\in \mathcal{P}_m(\C, \C^{s\times t} \rbrace
\end{equation}
where $\set P_m(\C, \C^{s\times t}$ is the space of matrix polynomials of degree $m$.
If we then seek to perform a matrix-valued rational approximation using the Frobenius norm
\begin{equation}
	\minimize_{\ma R \in \set R_{m,n}(\C,\C^{s\times t})}
		\sum_{j=1}^N \| \ma R(z_j) - \ma F(z_j) \|_\fro
\end{equation} 
we can easily generalize our approaches to this problem.
For example, considering a pole-residue expansion of a degree $(n-1,n)$ rational function,
$\ma R(z)$,
\begin{equation}
	\ma R(z) = \sum_{k=1}^n \frac{ \ve \rho_k}{z - \lambda_k},
	\quad \ma R: \C \to \C^{s\times t}, \ve \rho_k \in \C^{s\times t}.
\end{equation}
we can similarly use VARPRO to implicitly solve for linear parameters $\ve \rho_k$.
However, in the context of model reduction,
it is often desired to impose a constraint that each $\ve \rho_k$
is rank-1, so that the dimension of the reduced order model is $n$.
It is an open question as to how best to incorporate this constraint.

\section*{Acknowledgements}
The first author would like to thank Paul Constantine for his support 
during the writing of this manuscript.


\bibliographystyle{siamplain}
\bibliography{ratfit}
\end{document}